\documentclass[11pt]{amsart}   

\newtheorem{theorem}{Theorem}[section]
\newtheorem{lemma}[theorem]{Lemma}
\newtheorem{proposition}[theorem]{Proposition}
\theoremstyle{definition}
\newtheorem{definition}[theorem]{Definition}
\theoremstyle{remark}
\newtheorem{remark}[theorem]{Remark}
\newtheorem{problem}[theorem]{Problem}
\numberwithin{equation}{section}

\def\ep{\varepsilon}

\def\lip{\hskip0.02cm{\rm Lip}\hskip0.01cm}

\def\span{\hskip0.02cm{\rm span}\hskip0.01cm}
\def\mju{\mathcal{U}}
\def\mv{\mathcal{V}}

\begin{document}

\title{Coarse embeddability into Banach spaces}

\author{M.I.~Ostrovskii}
\address{Department of Mathematics and Computer Science\\
St. John's University\\
8000 Utopia Parkway\\
Queens, NY 11439, USA}
\email{ostrovsm@stjohns.edu}




\begin{abstract} The main purposes of this paper are (1) To survey the area of
coarse embeddability of metric spaces into Banach spaces, and, in
particular, coarse embeddability of different Banach spaces into
each other; (2) To present new results on the problems: (a)
Whether coarse non-embeddability into $\ell_2$ implies presence of
expander-like structures? (b) To what extent $\ell_2$ is the most
difficult space to embed into?
\end{abstract}

\maketitle

\begin{large}

\section{Introduction}

\subsection{Basic definitions}

Let $A$ and $B$ be metric spaces with metrics $d_A$ and $d_B$,
respectively.

\begin{definition}\label{coarse} {\rm A mapping
$f:A\to B$ is called a {\it coarse embedding} (or a {\it uniform
embedding}) if there exist functions
$\rho_1,\rho_2:[0,\infty)\to[0,\infty)$ such that

\noindent {\bf 1.} $\forall x,y\in A~ \rho_1(d_A(x,y))\le
d_B(f(x),f(y))\le\rho_2(d_A(x,y))$.

\noindent {\bf 2.} $\lim_{r\to\infty}\rho_1(r)=\infty$.}
\end{definition}

\begin{remark}
We prefer to use the term {\it coarse embedding} because in the
Nonlinear Functional Analysis the term {\it uniform embedding} is
used for uniformly continuous injective maps who\-se inverses are
uniformly continuous on their domains of definition, see
\cite[p.~3]{BL00}. In some of the papers cited below the term {\it
uniform embedding} is used.
\end{remark}

\begin{definition} {\rm A mapping $f:A\to B$ is called
{\it Lipschitz} if there exists a constant $0\le L<\infty$
such that
\begin{equation}\label{Lip}d_B(f(x),f(y))\le L\cdot d_A(x,y).\end{equation} The infimum of all $L>0$ for which the inequality in (\ref{Lip}) is valid is called the {\it Lipschitz constant} of $f$ and is denoted by $\lip(f)$.
A Lipschitz mapping is called a {\it Lipschitz embedding} if it is
one-to-one, and its inverse, defined as a mapping from the image
of $f$ into $A$, is also a Lipschitz mapping.}
\end{definition}

\begin{definition}\label{BG} {\rm A metric space $A$ is said to have {\it bounded geometry} if for
each $r>0$ there exist a positive integer $M(r)$ such that each
ball $B(x,r)=\{y\in A:~d_A(x,y)\le r\}$ of radius $r$ contains at
most $M(r)$ elements.}
\end{definition}

\begin{definition}\label{D:LF} {\rm A metric space is called
{\it locally finite} if all balls in it have finitely many
elements.}
\end{definition}

Our terminology and notation of Banach space theory follows
\cite{BL00} and \cite{JL01}.

\subsection{Some history and motivation}

M.~Gromov \cite{Gro95} suggested to use coarse embeddings of
Cayley graphs of infinite groups with finitely many generators and
finitely many relations (with their graph-theoretical metric) into
a Hilbert space or into a uniformly convex Banach space as a tool
for working on such well-known conjectures as the Novikov
conjecture and the Baum--Connes conjecture (discussion of these
conjectures is beyond the scope of this paper). G.~Yu \cite{Yu00}
and G.~Kasparov and G.~Yu \cite{KY06} have shown that this is
indeed a very powerful tool. G.~Yu \cite{Yu00} used the condition
of coarse embeddability of metric spaces with bounded geometry
into a Hilbert space; G.~Kasparov and G.~Yu \cite{KY06} used the
condition of coarse embeddability of metric spaces with bounded
geometry into a uniformly convex Banach space. These results made
the following problem posed by M.~Gromov in \cite[Problem
(4)]{Gro95} very important:
\smallskip

{\it ``Does every finitely generated or finitely presented group
admit a uniformly metrically proper Lipschitz embedding into a
Hilbert space? Even such an embedding into a reflexive uniformly
convex Banach space would be interesting. This seems hard.''}
\medskip

Also, they attracted attention to the following generalized
version of the problem:
\smallskip

{\it Whether each metric space with bounded geometry is coarsely
embeddable into a uniformly convex Banach space? }
\medskip

The result of G.~Kasparov and G.~Yu \cite{KY06} also made it
interesting to compare classes of metric spaces embeddable into
different Banach spaces (with particular interest to spaces with
bounded geometry).

\section{Obstructions to embeddability of spaces with bounded
geometry}

M.~Gromov \cite[Remark (b), p.~218]{Gro93} wrote: ``There is no
known geometric obstruction for uniform embeddings into infinite
dimensional Banach spaces.'' Writing this M.~Gromov was unaware of
P.~Enflo's work \cite{Enf69} in which it was shown that there is
no uniformly continuous embedding with uniformly continuous
inverse of the Banach space $c_0$ into a Hilbert space.
A.N.~Dranishnikov, G. Gong, V.~Lafforgue, and G.~Yu \cite[Section
6]{DGLY02} adjusted the construction of P.~Enflo \cite{Enf69} in
order to prove that there exist locally finite metric spaces which
are not coarsely embeddable into Hilbert spaces. After
\cite{DGLY02} was written, M.~Gromov (see \cite[p.~158]{Gro00})
observed that expanders provide examples of spaces with bounded
geometry which are not coarsely embeddable into a Hilbert space
and into $\ell_p$ for $1\le p<\infty$. Recall the definition (see
\cite{DSV03} for an accessible introduction to the theory of
expanders).

\begin{definition} {\rm For a finite graph $G$ with vertex set $V$ and a subset $F\subset V$
by $\partial F$ we denote the set of edges connecting $F$ and
$V\backslash F$. The {\it expanding constant} (also known as {\it
Cheeger constant}) of $G$ is
$$h(G)=\inf\left\{\frac{|\partial F|}{|F|}:~ F\subset V,~ 0<|F|\le |V|/2\right\}.$$

A sequence $\{G_n\}$ of graphs is called a {\it family of
expanders} if all of $G_n$ are finite, connected, $k$-regular for
some $k\in\mathbb{N}$ (that is, each vertex is adjacent to exactly
$k$ other vertices), their expanding constants $h(G_n)$ are
bounded away from $0$ (that is, there exists $\ep>0$ such that
$h(G_n)\ge\ep$ for all $n$), and their orders (numbers of
vertices) tend to $\infty$ as $n\to\infty$.}
\end{definition}

We consider (vertex sets of) connected graphs as metric spaces,
with their standard graph-theoretic distance: the distance between
two vertices is the number of edges in the shortest path joining
them.
\medskip

Let $A$ be a metric space containing isometric copies of all
graphs from some family of expanders. The Gromov's observation is:
$A$ does not embed coarsely into $\ell_p$ for $1\le p<\infty$
(see \cite[pp.~160--161]{Roe03} for a detailed proof, it is worth
mentioning that the result can be proved using the argument which
is well-known in the theory of Lipschitz embeddings of finite
metric spaces, see \cite[pp.~192--193]{Mat97}).
\medskip

M.~Gromov \cite{Gro00} suggested to use random groups in order to
prove that there exist Cayley graphs of finitely presented groups
which are not coarsely embeddable into a Hilbert space. Many
details on this approach were given in the paper M.~Gromov
\cite{Gro03} (some details were explained in \cite{Ghy04},
\cite{Oll05}, and \cite{Sil03}). However, to the best of my
knowledge, the work on clarification of all of the details of the
M.~Gromov's construction has not been completed (as of now).
\medskip

The posed above problem about the existence of coarse embeddings
of spaces with bounded geometry into uniformly convex Banach
spaces was recently solved in the negative by V.~La\-fforgue
\cite{Laf07+}, his construction is also expander-based.
\medskip

N.~Ozawa \cite[Theorem A.1]{Oza04} proved that a metric space $A$
containing isometric copies of all graphs from some family of
expanders does not embed coarsely into any Banach space $X$ such
that $B_X$ (the unit ball of $X$) is uniformly homeomorphic to a
subset of a Hilbert space. See \cite[Chapter 9, Section 2]{BL00}
for results on spaces $X$ such that $B_X$ is uniformly homeomorphic to $B_{\ell_2}$.
\medskip

It would be very interesting to find out whether each metric space
with bounded geometry which is not coarsely embeddable into a
Hilbert space contains a substructure similar to a family of
expanders. A version of this problem was posed in \cite{GK04}
using the following terminology:

\begin{definition} {\rm A metric space $X$ {\it weakly contains} a
family $\{G_n\}_{n=1}^\infty$ of expanders with vertex sets
$\{V_n\}_{n=1}^\infty$ if there are maps $f_n:V_n\to Y$ satisfying
\begin{itemize}
\item[(i)] $\displaystyle{\sup_{n}\lip(f_n)<\infty}$, \item[(ii)]
$\displaystyle{\lim_{n\to\infty}\sup_{v\in
V_n}}\frac{|f_n^{-1}(v)|}{|V_n|}=0$.
\end{itemize}
}\end{definition}

\begin{problem} \cite[p.~261]{GK04} Let $\Gamma$ be a
finitely presented group whose Cayley graph $G$ with its natural
metric is not coarsely embeddable into $\ell_2$. Does it follow
that $G$ weakly contains a family of expanders?
\end{problem}

The following theorem can help to finding an expander-like
structure in metric spaces with bounded geometry which are not
coarsely embeddable into a Hilbert space. In the theorem we
consider coarse embeddability into $L_1=L_1(0,1)$. For technical
reasons it is more convenient to work with $L_1$. As we shall see
in section \ref{S:class} coarse embeddability into $L_1$ is
equivalent to coarse embeddability into a Hilbert space.

\begin{theorem}\label{T:mu} Let $M$ be a locally finite metric space which is not coarsely embeddable into $L_1$. Then
there exists a constant $D$, depending on $M$ only, such that for
each $n\in\mathbb{N}$ there exists a finite set $B_n\subset
M\times M$ and a probability measure $\mu$ on $B_n$ such that
\begin{itemize}
\item $d_M(u,v)\ge n$ for each $(u,v)\in B_n$.

\item For each Lipschitz function $f:M\to L_1$ the inequality
\begin{equation}\label{Obs5}\int_{B_n}||f(u)-f(v)||_{L_1}d\mu(u,v)\le D\lip(f)\end{equation}
holds.
\end{itemize}
\end{theorem}

\begin{lemma}\label{L:1} There exists a constant $C$ depending on $M$ only
such that for each Lipschitz function $f:M\to L_1$ there exists a
subset $B_f\subset M\times M$ such that
$\displaystyle{\sup_{(x,y)\in B_f}d_M(x,y)}=\infty$, but
$\displaystyle{\sup_{(x,y)\in B_f}||f(x)-f(y)||_{L_1}}\le
C\lip(f)$.
\end{lemma}

\begin{proof} Assume the contrary. Then,
for each $n\in\mathbb{N}$, the number $n^3$ cannot serve as $C$.
This means, that for each $n\in\mathbb{N}$ there exists a
Lipschitz mapping $f_n:M\to L_1$ such that for each subset
$U\subset M\times M$ with
$$\sup_{(x,y)\in U}d_M(x,y)=\infty,$$
we have
$$\sup_{(x,y)\in U}||f_n(x)-f_n(y)||> n^3\lip(f_n).$$
We choose a point in $M$ and denote it by $O$. Without loss of
generality we may assume that $f_n(O)=0$. Consider the mapping
$$f:M\to\left(\sum_{n=1}^\infty\oplus L_1\right)_1\subset L_1$$
given by
$$f(x)=\sum_{n=1}^\infty\frac1{Kn^2}\cdot\frac{f_n(x)}{\lip(f_n)},$$
where $K=\sum_{n=1}^\infty\frac1{n^2}$. It is clear that the
series converges and $\lip(f)\le1$.
\medskip

Let us show that $f$ is a coarse embedding. We need an estimate
from below only (the estimate from above is satisfied because $f$ is
Lipschitz).
\medskip

The assumption implies that for each $n\in\mathbb{N}$ there is
$N\in\mathbb{N}$ such that
$$d_M(x,y)\ge N\Rightarrow ||f_n(x)-f_n(y)||>n^3\lip(f_n).$$
On the other hand
$$||f_n(x)-f_n(y)||> n^3\lip(f_n)\Rightarrow ||f(x)-f(y)||>\frac{n}{K}$$
Hence $f:M\to L_1$ is a coarse embedding and we get a
contradiction.
\end{proof}

\begin{lemma}\label{L:2} Let $C$ be the constant whose existence is proved in Lemma {\rm \ref{L:1}} and let $\ep$
be an arbitrary positive number. For each $n\in\mathbb{N}$ we can
find a finite subset $M_n\subset M$ such that for each Lipschitz
mapping $f:M\to L_1$ there is a pair $(u_{f,n},v_{f,n})\in
M_n\times M_n$ such that
\begin{itemize}
\item $d_M(u_{f,n},v_{f,n})\ge n$. \item
$||f(u_{f,n})-f(v_{f,n})||\le (C+\ep)\lip(f)$.
\end{itemize}
\end{lemma}

\begin{proof} We choose a point in $M$ and denote it by $O$. The ball in $M$ of radius $R$ centered at $O$
will be denoted by $B(R)$. It is clear that it suffices to prove
the result for $1$-Lipschitz mappings satisfying $f(O)=0$.
\medskip

Assume the contrary. Since $M$ is locally finite, this implies
that for each $R\in\mathbb{N}$ there is a 1-Lipschitz mapping
$f_R:M\to L_1$ such that $f_R(O)=0$ and, for $u,v\in B(R)$, the
inequality $d_M(u,v)\ge n$ implies
$||f_R(u)-f_R(v)||_{L_1}>C+\ep$.
\medskip

We refer to \cite{DK72}, \cite{Hei80}, or \cite[Chapter 8]{DJT95}
for results on ultraproducts, our terminology and notation follows
\cite{DJT95}. We form an ultraproduct of the mappings
$\{f_R\}_{R=1}^\infty$, that is, a mapping $f:M\to
(L_1)^{\mathcal{U}}$, given by $f(m)=\{f_R(m)\}_{R=1}^\infty$,
where $\mathcal{U}$ is a non-trivial ultrafilter on $\mathbb{N}$
and $(L_1)^{\mathcal{U}}$ is the corresponding ultrapower.  Each
ultrapower of $L_1$ is isometric to an $L_1$ space on some measure
space (see \cite[Theorem 8.7]{DJT95}, \cite{DK72}, \cite{Hei80}),
and its separable subspaces are isometric to subspaces of
$L_1(0,1)$ (see \cite[p.~168]{DS58}, \cite[pp.~14--15]{JL01}, and
\cite[pp.~399 \&\ 416]{Roy88}). Therefore we can consider $f$ as a
mapping into $L_1(0,1)$. It is easy to verify that $\lip(f)\le 1$
and that $f$ satisfies the condition
$$d_M(u,v)\ge n\Rightarrow ||f(u)-f(v)||_{L_1}\ge (C+\ep).$$
We get a contradiction with the definition of $C$.
\end{proof}

\begin{proof}[Proof of Theorem {\rm \ref{T:mu}}] Let $D$ be a
number satisfying $D>C$, and let $B$ be a number satisfying
$C<B<D$.

According to Lemma \ref{L:2}, there is a finite subset $M_n\subset
M$ such that for each $1$-Lipschitz function $f$ on $M$ there is a
pair $(u,v)$ in $M_n$ such that $d_M(u,v)\ge n$ and
$||f(u)-f(v)||\le B$.
\medskip

Let $\alpha_n$ be the cardinality of $M_n$, we choose a point in
$M_n$ and denote it by $O$. Proving the theorem it is enough to
consider $1$-Lipschitz functions $f:M_n\to L_1$ satisfying
$f(O)=0$. Each $\alpha_n$-element subset of $L_1$ is isometric to
a subset in $\ell_1^{\alpha_n(\alpha_n-1)/2}$ (see \cite{Wit86},
\cite{Bal90}). Therefore it suffices to prove the result for
$1$-Lipschitz embeddings into $\ell_1^{\alpha_n(\alpha_n-1)/2}$.
It is clear that it suffices to prove the inequality
\[\int_{B_n}||f(u)-f(v)||d\mu(u,v)\le
B\] for a $\left(\frac{D-B}2\right)$-net in the set of all functions
satisfying the conditions mentioned above, endowed with the metric
$$\tau(f,g)=\max_{m\in M_n}||f(m)-g(m)||$$

By compactness there exists a finite net satisfying the condition.
Let $N$ be such a net. We are going to use the minimax theorem,
see, e.g. \cite[p.~344]{Str80}. In particular, we use the notation
similar to the one used in \cite{Str80}.
\medskip

Let $A$ be the matrix whose columns are labelled by functions from
$N$, whose rows are labelled by pairs $(u,v)$ of elements of $M_n$
satisfying $d_M(u,v)\ge n$, and whose entry on the intersection of
the column corresponding to $f$, and the row corresponding to
$(u,v)$ is $||f(u)-f(v)||$.
\medskip

Then, for each column vector $x=\{x_f\}_{f\in N}$ with $x_f\ge 0$
and $\sum_{f\in N} x_f=1$, the entries of the product $Ax$ are the
differences $||F(u)-F(v)||$, where
$F:M\to\displaystyle{\left(\sum_{f\in N}\oplus
\ell_1^{\alpha_n(\alpha_n-1)/2}\right)_1}$ is given by
$F(m)=\displaystyle{\sum_{f\in N} x_ff(m)}$. The function $F$ can
be considered as a function into $L_1$. It satisfies $\lip(F)\le
1$. Hence there is a pair $(u,v)$ in $M_n$ satisfying $d_M(u,v)\ge
n$ and $||F(u)-F(v)||\le B$. Therefore we have
$$\max_x\min_\mu\mu Ax\le B,$$
where the minimum is taken over all vectors $\mu=\{\mu(u,v)\}$,
indexed by $u,v\in M_n$, $d_M(u,v)\ge n$, and satisfying the
conditions $\mu(u,v)\ge 0$ and $\sum \mu(u,v)=1$.
\medskip

By the von Neumann minimax theorem \cite[p.~344]{Str80}, we have
$$\min_\mu\max_x\mu Ax\le B,$$
which is exactly the inequality we need to prove because $\mu$ can
be regarded as a probability measure on the set of pairs from
$M_n$ with distance $\ge n$.
\end{proof}

\section{Coarse embeddability into reflexive Banach spaces}

The first result of this nature was obtained by N.~Brown and
E.~Guentner \cite[Theorem 1]{BG05}. They proved that for each
metric space $A$ having bounded geometry there is a sequence
$\{p_n\}$, $p_n>1$, $\lim_{n\to\infty} p_n=\infty$ such that $A$
embeds coarsely into the Banach space
$(\sum_{n=1}^\infty\oplus\ell_{p_n})_2$, which is, obviously,
reflexive.
\medskip

This result was strengthened in \cite{BL07+}, \cite{Kal07+}, and
\cite{Ost06a}. (Observe that the space
$(\sum_{n=1}^\infty\oplus\ell_{p_n})_2$ has no cotype.)

\begin{theorem}{\rm \cite{Ost06a}}
Let $X$ be a Banach space with no cotype and let $A$ be a locally
finite metric space. Then $A$ embeds coarsely into $X$.
\end{theorem}

\begin{theorem}{\rm \cite{BL07+}}
Let $X$ be a Banach space with no cotype and let $A$ be a locally
finite metric space. Then there exists a Lipschitz embedding of
$A$ into $X$.
\end{theorem}

\begin{remark}
Interested readers can reconstruct the proof from \cite{BL07+} by
applying Proposition \ref{Z} (see below) to $Z=c_0$ in combination
with the result of I.~Aharoni mentioned in Section \ref{S:CE}.
\end{remark}

\begin{definition} {\rm A metric space $(X,d)$ is called
{\it stable} if for any two bounded sequences $\{x_n\}$ and
$\{y_n\}$ in $X$ and for any two non-trivial ultrafilters $\mju$
and $\mv$ on $\mathbb{N}$ the condition
$$\displaystyle{\lim_{n,\mju}\lim_{m,\mv}d(x_n,y_m)=\lim_{m,\mv}\lim_{n,\mju}d(x_n,y_m)}$$
holds.}
\end{definition}

\begin{theorem}{\rm \cite{Kal07+}} Let $A$ be a stable metric space. Then $A$ embeds coarsely
into a reflexive Banach space.
\end{theorem}

\begin{remark} It is easy to see that locally finite metric spaces are stable.
\end{remark}

N.J.~Kalton \cite{Kal07+} found examples of Banach spaces which
are not coarsely embeddable into reflexive Banach spaces, $c_0$ is
one of the examples of such spaces. Apparently his result provides
the first example of a metric spaces which is not coarsely
embeddable into reflexive Banach spaces. (See the Problem (3) in
the list of open problems in \cite{Pes07}.)

\section{Coarse classification of Banach spaces}\label{S:class}

As we already mentioned the result of G.~Kasparov and G.~Yu
\cite{KY06} makes it very interesting to compare the conditions of
coarse embeddability into a Banach space $X$ for different spaces
$X$. Since compositions of coarse embeddings are coarse
embeddings, one can approach this problem by studying coarse
embeddability of Banach space into each other. In this subsection
we describe the existing knowledge on this matter.

\subsection{Essentially nonlinear coarse embeddings}\label{S:CE}

There are many examples of pairs $(X,Y)$ of Banach spaces such
that $X$ is coarsely embeddable into $Y$, but the
Banach-space-theoretical structure of $X$ is quite different from
the Banach-space-theore\-ti\-cal structure of each subspace of
$Y$:
\begin{itemize}
\item A result which goes back to I.J.~Schoenberg \cite{Sch38}
(see \cite[p.~385]{Mat07+} for a simple proof) states that $L_1$
with the metric $\sqrt{||x-y||_1}$ is isometric to a subset of
$L_2$. Hence $L_1$ and all of its subspaces, in particular, $L_p$
and $\ell_p~ (1\le p\le 2)$ (see \cite{Kad58} and \cite{BDK66})
embed coarsely into $L_2=\ell_2$.

\item This result was generalized by M.~Mendel and A.~Naor
\cite[Remark 5.10]{MN04}: For every $1\le q<p$ the metric space
$(L_q, ||x-y||_{L_q}^{q/p})$ is isometric to a subspace of $L_p$.

\item The well-known result of I.~Aharoni \cite{Aha74} implies
that  each separable Banach space is coarsely embeddable into
$c_0$ (although its Banach space theoretical properties can be
quite different from those of any subspace of $c_0$). A simpler
proof of this result was obtained in \cite{Ass78}, see, also,
\cite[p.~176]{BL00}.

\item N.J.~Kalton \cite{Kal04} proved that $c_0$ embeds coarsely
into a Banach space with the Schur property.

\item P.~Nowak \cite{Now06} proved that $\ell_2$ is coarsely
embeddable into $\ell_p$ for all $1\le p\le \infty$.
\end{itemize}

\subsection{Obstructions to coarse embeddability of Banach spaces}\label{S:???}

The list of discovered obstructions to coarse
embeddability also constantly increases:
\begin{itemize}

\item Only minor adjustments of the argument of Y.~Raynaud
\cite{Ray83} (see, also \cite[pp.~212--215]{BL00}) are needed to
prove the following results:

\begin{itemize}

\item[(1)] Let $A$ be a Banach space with a spreading basis which
is not an unconditional basis. Then $A$ does not embed coarsely
into a stable metric space. (See \cite[p.~429]{BL00} for the
definition of a spreading basis and \cite{KM81} for examples of
stable Banach spaces. Examples of stable Banach spaces include
$L_p$ $(1\le p<\infty)$.)

\item[(2)] Let $A$ be a nonreflexive Banach space with non-trivial
type. Then $A$ does not embed coarsely into a stable metric space.
(Examples of nonreflexive Banach spaces with non-trivial type were
constructed in \cite{Jam74}, \cite{JL75}, \cite{PX87}.)

\end{itemize}
\item A.N.~Dranishnikov, G.~Gong, V.~Lafforgue, and G.~Yu
\cite{DGLY02} adjusted the argument of P.~Enflo \cite{Enf69} to
prove that Banach spaces with no cotype are not coarsely
embeddable into $\ell_2$.

\item W.~B.~Johnson and L.~Randrianarivony \cite{JR06} proved that
$\ell_p$ $(p>2)$ is not coarsely embeddable into $\ell_2$.

\item M.~Mendel and A.~Naor \cite{MN07+} proved (for $K$-convex
spa\-ces) that cotype of a Banach space is an obstruction to
coarse embeddability, in particular, $\ell_p$ is not coarsely
embeddable into $\ell_q$ when $p>q\ge 2$.

\item L.~Randrianarivony \cite{Ran06} strengthened the result from
\cite{JR06} to a characterization of quasi-Banach spaces which
embed coarsely into a Hilbert space, and proved: a separable
Banach space is coarsely embeddable into a Hilbert space if and
only if it is isomorphic to a subspace of $L_0(\mu)$.

\item N.J.~Kalton \cite{Kal07+} found some more obstructions to
coarse embeddability. In particular, N.J.~Kalton discovered an
invariant, which he named the $\mathcal{Q}$-property, which is
necessary for coarse embeddability into reflexive Banach spaces.
\end{itemize}

\subsection{To what extent is $\ell_2$ the most difficult space to
embed into?}

Because $\ell_2$ is, in many respects, the `best' space, and
because of Dvoretzky's theorem (see \cite{Dvo61} and \cite{MS86})
it is natural to expect that $\ell_2$ is among the most difficult
spaces to embed into. The strongest possible result in this
direction would be a positive solution to the following problem.

\begin{problem} Does $\ell_2$ embed coarsely into an arbitrary
infinite dimensional Banach spa\-ce?
\end{problem}

This problem is still open, but the coarse embeddability of
$\ell_2$ is known for wide classes of Banach spaces. As was
mentioned above, P.W.~Nowak \cite{Now06} proved that $\ell_2$
embeds coarsely into $\ell_p$ for each $1\le p\le\infty$. In
Section \ref{L2} we prove that $\ell_2$ embeds coarsely into a
Banach space containing a subspace with an unconditional basis
which does not contain $\ell_\infty^n$ uniformly (Theorem
\ref{T:L2}). This result is a generalization of P.W.~Nowak's
result mentioned above because the spaces $\ell_p$ $(1\le
p<\infty)$ satisfy the condition of Theorem \ref{T:L2}, but the
spaces satisfying the condition of Theorem \ref{T:L2} do not
necessarily contain subspaces isomorphic to $\ell_p$ (see
\cite{FT74}, and \cite[Section 2.e]{LT77}).
\medskip

In all existing applications of coarse embeddability results the
most important is the case when we embed spaces with bounded
geometry into Banach spaces. In this connection the following
result from \cite{Ost06b} is of interest.

\begin{theorem}[\cite{Ost06b}] Let $A$ be a locally finite metric space which embeds coarsely into a
Hilbert space, and let $X$ be an infinite dimensional Banach
space. Then there exists a coarse embedding $f:A\to X$.
\end{theorem}

In this paper we use an idea of F.~Baudier and G.~Lancien
\cite{BL07+}, and prove this result in a stronger form, for
Lipschitz embeddings (see Section \ref{S:LF}):

\begin{theorem}\label{T:HLF} Let $M$ be a locally finite subset of
a Hilbert space. Then $M$ is Lipschitz embeddable into an
arbitrary infinite dimensional Banach space.
\end{theorem}

\section{Coarse embeddings of $\ell_2$}\label{L2}

\begin{theorem}\label{T:L2} Let $X$ be a Banach space containing a
subspace with an unconditional basis which does not contain
$\ell_\infty^n$ uniformly. Then $\ell_2$ embeds coarsely into $X$.
\end{theorem}

\begin{proof} We use the criterion for coarse embeddability into a Hilbert space due to
M.~Dadarlat and E.~Guentner \cite[Proposition 2.1]{DG03} (see
\cite{LW06} and \cite{Now06} for related results). We state it as
a lemma (by $S(X)$ we denote the unit sphere of a Banach space
$X$).

\begin{lemma}[\cite{DG03}]\label{DG} A metric space $A$ admits a coarse
embedding into $\ell_2$ if and only if for every $\ep>0$ and every
$R>0$ there exists a map $\zeta:A\to S(\ell_2)$ such that
\begin{itemize}
\item[{\rm (i)}] $d_A(x,y)\le R$ implies
$||\zeta(x)-\zeta(y)||\le\ep$.

\item[{\rm (ii)}] $\lim_{t\to\infty}\inf\{||\zeta(x)-\zeta(y)||:~
x,y\in A,~ d_A(x,y)\ge t\}=\sqrt{2}$.
\end{itemize}
\end{lemma}

We assume without loss of generality that $X$ has an unconditional
basis $\{e_i\}_{i\in\mathbb{N}}$. Let
$\mathbb{N}=\cup_{i=1}^\infty \mathbb{N}_i$ be a partition of
$\mathbb{N}$ into infinitely many infinite subsets. Let
$X_i=\hbox{cl}(\span\{e_i\}_{i\in\mathbb{N}_i})$. By the theorem
of E.~Odell and T.~Schlumprecht \cite{OS94} (see, also,
\cite[Theorem 9.4]{BL00}), for each $i\in \mathbb{N}$ there exists
a uniform homeomorphism $\varphi_i:S(\ell_2)\to S(X_i)$. We apply
Lemma \ref{DG} in the case when $A=\ell_2$. By the uniform
continuity of $\varphi_i$ and $\varphi_i^{-1}$ we get: for each
$i\in \mathbb{N}$ there exists $\delta_i>0$ and a map
$\zeta_i:\ell_2\to S(X_i)$ such that

\begin{equation}\label{E:ii}
\lim_{t\to\infty}\inf\{||\zeta_i(x)-\zeta_i(y)||_{X_i}:~
||x-y||_{\ell_2}\ge t\}\ge\delta_i.
\end{equation}

\begin{equation}\label{E:i}
||x-y||_{\ell_2}\le i\hbox{ implies }
||\zeta_i(x)-\zeta_i(y)||_{X_i}\le\frac{\delta_i}{i2^i}.
\end{equation}

Fix $x_0\in\ell_2$. Let $f:\ell_2\to X$ be the map defined as the
direct sum of the maps
$\displaystyle{\frac{i}{\delta_i}\left(\zeta_i(x)-\zeta_i(x_0)\right)}$.
We claim that it is a coarse embedding (the fact that it is a
well-defined map follows from (\ref{E:i})).
\medskip

Let $||x-y||=r$, then for $i\ge r$ we get
$||\frac{i}{\delta_i}\zeta_i(x)-\frac{i}{\delta_i}\zeta_i(y)||_{X_i}\le
\frac1{2^i}$. Hence $||f(x)-f(y)||\le\sum_{i=1}^{\lceil
r\rceil-1}\frac{2i}{\delta_i}+\sum_{i=\lceil
r\rceil}^\infty\frac1{2^i}=:\rho_2(r)$. We proved an estimate from
above.
\medskip

To prove an estimate from below, it is enough, for a given
$h\in\mathbb{R}$, to find $t\in\mathbb{R}$ such that
$||x-y||_{\ell_2}\ge t$ implies $||f(x)-f(y)||_X\ge h$. For this,
by unconditionality (we assume, for simplicity, that the basis of
$X$ is $1$-unconditional), it is enough to find $i\in\mathbb{N}$
such that $||x-y||_{\ell_2}\ge t$ implies
$||\frac{i}{\delta_i}\zeta_i(x)-\frac{i}{\delta_i}\zeta_i(y)||_{X_i}\ge
h$. We choose an arbitrary $i>h$. The conclusion follows from the
condition (\ref{E:ii}).
\end{proof}

\section{Lipschitz embeddings of locally finite metric spaces}\label{S:LF}

The purpose of this section is to prove Theorem \ref{T:HLF}. We
prove the main step in our argument (Proposition \ref{Z}) in a
somewhat more general context than is needed for Theorem
\ref{T:HLF}, because it can be applied in some other situations
(see, in this connection, the paper \cite{Bau07+} containing two
versions of Proposition \ref{Z}). The coarse version of this
result was proved in \cite{Ost06b}, in the proof of the Lipschitz
version we use an idea from \cite{BL07+}.

\begin{proposition}\label{Z} Let $A$ be a locally finite subset of a
Banach space $Z$. Then there exists a sequence of finite
dimensional linear subspaces $Z_i~(i\in\mathbb{N})$ of $Z$ such
that $A$ is Lipschitz embeddable into each Banach space $Y$ having
a finite dimensional Schauder decomposition $\{Y_i\}_{i=1}^\infty$
with $Y_i$ linearly isometric to $Z_i$.
\end{proposition}

See \cite[Section 1.g]{LT77} for information on Schauder
decompositions. It is clear that we may restrict ourselves to the
case when the Schauder decomposition satisfies
\begin{equation}\label{mon}
||y_i||\le\left\|\sum_{i=1}^\infty y_i\right\|\hbox{ when }y_i\in
Y_i~\hbox{ for each }~ i\in\mathbb{N}.
\end{equation}

\begin{proof}
Let $Z_i$ be the linear subspace of $Z$ spanned by $\{a\in
A:~||a||_Z\le 2^i\}$ and let $S_i=\{a\in A:~2^{i-1}\le ||a||_Z\le
2^i\}$. Let $T_i:Z_i\to Y_i$ be some linear isometries and let
$E_i:Z_i\to Y$ be compositions of these linear isometries with the
natural embeddings $Y_i\to Y$. We define an embedding
$\varphi:A\to Y$ by
$$\varphi(a)=\frac{2^i-||a||_Z}{2^{i-1}}E_i(a)+\frac{||a||_Z-2^{i-1}}{2^{i-1}}E_{i+1}(a)\hbox{ for }a\in S_i.$$
One can check that there is no ambiguity for $||a||_Z=2^i$.

\begin{remark} The mapping $\varphi$ is a straightforward generalization of the
mapping constructed in \cite{BL07+}.
\end{remark}

It remains to verify that $\varphi$ is a Lipschitz embedding. We
consider three cases.

\begin{itemize}
\item[(1)] $a,b$ are in the same $S_i$; \item[(2)] $a,b$ are in
consecutive sets $S_i$, that is, $b\in S_i$, $a\in S_{i+1}$;
\item[(3)] $a,b$ are in `distant' sets $S_i$, that is, $b\in S_i$,
$a\in S_k$, $k\ge i+2$.
\end{itemize}

Everywhere in the proof we assume $||a||\ge ||b||$.
\medskip

\noindent{\sc Case (1).} The inequality (\ref{mon}) implies that
the number $$||\varphi(a)-\varphi(b)||_Y$$ is between the maximum
and the sum of the numbers
\begin{equation}\label{i}
\left\|\frac{2^i-||a||_Z}{2^{i-1}}E_i(a)-\frac{2^i-||b||_Z}{2^{i-1}}E_i(b)\right\|,
\end{equation}
\begin{equation}\label{i+1}
\left\|\frac{||a||_Z-2^{i-1}}{2^{i-1}}E_{i+1}(a)-\frac{||b||_Z-2^{i-1}}{2^{i-1}}E_{i+1}(b)
\right\|.
\end{equation}

It is clear that the norm in (\ref{i}) is between the numbers
\[\frac{2^i-||a||_Z}{2^{i-1}}||E_i(a)-E_i(b)||\mp\frac{||a||_Z-||b||_Z}{2^{i-1}}||E_{i+1}(b)||,\]
and the norm in (\ref{i+1}) is between the numbers
\[\frac{||a||_Z-2^{i-1}}{2^{i-1}}||E_{i+1}(a)-E_{i+1}(b)||\mp
\frac{||a||_Z-||b||_Z}{2^{i-1}}||E_{i+1}(b)||.\] Therefore
\[\begin{split}
\frac12\left(||a-b||_Z-\frac{||a||_Z-||b||_Z}{2^{i-2}}||b||_Z\right)\le||\varphi(a)-\varphi(b)||_Y\\
\le ||a-b||_Z+\frac{||a||_Z-||b||_Z}{2^{i-2}}||b||_Z.\end{split}\]

This inequality implies a suitable estimate from above for the
Lipschitz constant of $\varphi$, and an estimate for the Lipschitz
constant of its inverse in the case when $||a-b||_Z$ is much
larger than $||a||_Z-||b||_Z$, for example, if $||a-b||_Z\ge
5(||a||_Z-||b||_Z)$. To complete the proof in the case (1) it
suffices to estimate $||\varphi(a)-\varphi(b)||$ from below in the
case when $||a||_Z-||b||_Z\ge\displaystyle{\frac{||a-b||_Z}5}$. In
this case we use the observation that for $a,b\in S_i$ satisfying
$||a||_Z\ge ||b||_Z$ the sum of (\ref{i}) and (\ref{i+1}) can be
estimated form below by
\[\begin{split}
&\left(\frac{2^i-||a||_Z}{2^{i-1}}||a||_Z-\frac{2^i-||b||_Z}{2^{i-1}}||b||_Z\right)\\
+&{ }\left(\frac{||a||_Z-2^{i-1}}{2^{i-1}}||a||_Z-\frac{||b||_Z-2^{i-1}}{2^{i-1}}||b||_Z\right)\\
=&{~ }||a||_Z-||b||_Z\ge\frac{||a-b||_Z}5.\end{split}\] This
completes our proof in the case (1).
\medskip

\noindent{\sc Case (2).} The inequality (\ref{mon}) implies that
the number $$||\varphi(a)-\varphi(b)||_Y$$ is between the maximum
and the sum of the numbers
\begin{equation}\label{i2}
\left\|\frac{2^i-||b||_Z}{2^{i-1}}E_i(b)\right\|,
\end{equation}
\begin{equation}\label{i+12}
\left\|\frac{2^{i+1}-||a||_Z}{2^{i}}E_{i+1}(a)-\frac{||b||_Z-2^{i-1}}{2^{i-1}}E_{i+1}(b)\right\|,
\end{equation}
\begin{equation}\label{i+2}
\left\|\frac{||a||_Z-2^{i}}{2^{i}}E_{i+2}(a)\right\|.
\end{equation}

Both (\ref{i2}) and (\ref{i+2}) are estimated from above by
$2(||a||_Z-||b||_Z)$. As for (\ref{i+12}), we have \[
\left\|\frac{2^{i+1}-||a||_Z}{2^{i}}E_{i+1}(a)-\frac{||b||_Z-2^{i-1}}{2^{i-1}}E_{i+1}(b)\right\|\]
\begin{equation}\label{mid}=\left\|\frac{2^i-(||a||_Z-2^i)}{2^i}a+\frac{(2^i-||b||_Z)-2^{i-1}}{2^{i-1}}b
\right\|_Z\end{equation}
\[\le||a-b||_Z+2(||a||_Z-2^i)+2(2^i-||b||_Z)\le 3||a-b||_Z.
\]

We turn to estimate from below. From (\ref{i2}) and (\ref{i+2}) we
get
\[||\varphi(a)-\varphi(b)||\ge\max\{(2^i-||b||_Z),(||a||_Z-2^i)\}.\]
Therefore it suffices to find an estimate in the case when
\begin{equation}\label{sm}\max\{(2^i-||b||_Z),(||a||_Z-2^i)\}\le\frac{||a-b||_Z}5.\end{equation}
Rewriting (\ref{i+12}) in the same way as in (\ref{mid}), we get
\[||\varphi(a)-\varphi(b)||_Y\ge
\left\|(a-b)+\frac{2^i-||b||_Z}{2^{i-1}}b-\frac{||a||_Z-2^i}{2^i}a\right\|
\]
In the case when (\ref{sm}) is satisfied, we can continue this
chain of inequalities with
\[\ge||a-b||_Z-\frac45||a-b||_Z=\frac15||a-b||_Z.\]

\noindent{\sc Case (3).} In this case the number
$||\varphi(a)-\varphi(b)||_Y$ is between the maximum and the sum
of the four numbers:
$$\frac{2^i-||b||_Z}{2^{i-1}}||b||_Z,
\frac{||b||_Z-2^{i-1}}{2^{i-1}}||b||_Z,$$
$$
\frac{2^{k}-||a||_Z}{2^{k-1}}||a||_Z,
\frac{||a||_Z-2^{k-1}}{2^{k-1}}||a||_Z.
$$
Hence $||\varphi(a)-\varphi(b)||_Y$ is between
$\displaystyle{\frac{||a||_Z}2}$ (= the average of the last two
numbers) and $||a||_Z+||b||_Z$ (=the sum of all four numbers).

On the other hand,
$$\frac12||a||_Z\le||a||_Z-||b||_Z\le ||a-b||_Z\le||a||_Z+||b||_Z\le
2||a||_Z.$$ These inequalities immediately imply estimates for
Lipschitz constants.
\end{proof}

\begin{proof}[Proof of Theorem \ref{T:HLF}] Each
finite dimensional subspace of $\ell_2$ is isometric to $\ell_2^k$
for some $k\in\mathbb{N}$. By Proposition \ref{Z} there exists a
sequence $\{n_i\}_{i=1}^\infty$ such that $A$ embeds coarsely into
each Banach space $Y$ having a Schauder decomposition $\{Y_i\}$
with $Y_i$ isometric to $\ell_2^{n_i}$. On the other hand, using
Dvoretzky's theorem (\cite{Dvo61}, see, also, \cite[Section
5.8]{MS86}) and the standard techniques of constructing basic
sequences (see \cite[p.~4]{LT77}), it is easy to prove that for an
arbitrary sequence $\{n_i\}_{i=1}^\infty$ an arbitrary infinite
dimensional Banach space $X$ contains a subspace isomorphic to a
space having such Schauder decomposition.
\end{proof}

I would like to thank Tadeusz~Figiel and William B. Johnson for
useful conversations related to the subject of this paper.

\end{large}

\end{document}